\documentclass{gtmon_a}
\pdfoutput=1

\usepackage{pinlabel}
 
\proceedingstitle{The Zieschang Gedenkschrift}
\conferencestart{5 September 2007}
\conferenceend{8 September 2007}
\conferencename{Conference in honour of Heiner Zieschang}
\conferencelocation{Toulouse, France}

\editor{Michel Boileau}
\givenname{Michel}
\surname{Boileau}

\editor{Martin Scharlemann}
\givenname{Martin}
\surname{Scharlemann}

\editor{Richard Weidmann}
\givenname{Richard}
\surname{Weidmann}

\title[Configurations of torsion and Property $\rm{F}\mathbb R$]{A condition that prevents groups\\from acting nontrivially on trees}

\author[Martin R Bridson]{Martin R Bridson}
\givenname{Martin R}
\surname{Bridson}
\address{Mathematical Institute\\\newline
24-29 St Giles'\\
Oxford OX1 3LB\\UK}
\email{bridson@maths.ox.ac.uk}
\urladdr{http://www.maths.ox.ac.uk/~bridson}

\dedicatory{Dedicated to the memory of Heiner Zieschang}

\volumenumber{14}
\issuenumber{}
\publicationyear{2008}
\papernumber{08}
\startpage{129}
\endpage{133}

\doi{}
\MR{}
\Zbl{}

\arxivreference{} 

\keyword{property FA}
\keyword{automorphisms of free groups}
\keyword{$\mathbb R$--trees}
\subject{primary}{msc2000}{20E08}
\subject{primary}{msc2000}{20F65}

\received{20 October 2006}
\revised{26 January 2007}
\accepted{}
\published{29 April 2008}
\publishedonline{29 April 2008}
\proposed{}
\seconded{}
\corresponding{}
\version{}

\makeatletter
\def\cnewtheorem#1[#2]#3{\newtheorem{#1}{#3}[section]
\expandafter\let\csname c@#1\endcsname\c@thm}
\def\dnewtheorem#1[#2]#3{\newtheorem{#1}{#3}
\expandafter\let\csname c@#1\endcsname\c@thma}
\makeatother

\AtBeginDocument{\let\bar\wbar\let\tilde\wtilde\let\hat\what}
\AtBeginDocument{\def\S{\Sigma}}

\cnewtheorem{lemma}[thm]{Lemma}          
\cnewtheorem{coroll}[thm]{Corollary}
\newtheorem{thma}{Theorem}    
\dnewtheorem{lemmaa}[thma]{Lemma}          
\dnewtheorem{corolla}[thma]{Corollary}
\theoremstyle{definition}
\cnewtheorem{defn}[thm]{Definition}    
\cnewtheorem{remark}[thm]{Remark}
\newtheorem*{rem}{Remark}             
\makeautorefname{thma}{Theorem}
\makeautorefname{lemmaa}{Lemma}

\def\G{\Gamma}

\def\e{\varepsilon}
 
\def\An{\text{\rm{Aut}}(F_n)}
\def\Aut{\text{\rm{Aut}}}
\def\autn{\text{\rm{Aut}}(F_n)}

\def\slnz{\text{\rm{SL}}(n,\mathbb Z)}
\def\sln{\text{\rm{SL}}(n,\mathbb Z)}
\def\gln{\text{\rm{GL}}(n,\mathbb Z)}
\def\fix{\text{\rm{Fix}}}
\def\FR{\text{\rm{F}}\mathbb R}
\def\FA{\text{\rm{FA}}}
\def\B{\mathcal B}

\begin{document}

\begin{asciiabstract} We describe a simple criterion for showing
that a group has Serre's property FA. By exhibiting a
certain pattern of finite subgroups, we show that this criterion
is satisfied by
Aut(F_n) and SL(n,Z) when n>=3.
\end{asciiabstract}

\begin{htmlabstract}
We describe a simple criterion for showing
that a group has Serre's property FA.  By exhibiting a certain pattern
of finite subgroups, we show that this criterion is satisfied by
Aut(F<sub>n</sub>) and SL(n,<b>Z</b>) when n&ge; 3.
\end{htmlabstract}

\begin{abstract} We describe a simple criterion for showing
that a group has Serre's property FA.  By exhibiting a certain pattern
of finite subgroups, we show that this criterion is satisfied by
$\mathrm{Aut}(F_n)$ and $\mathrm{SL}(n,\mathbb{Z})$ when $n\ge 3$.
\end{abstract}

\maketitle


An $\mathbb R$--tree is a geodesic metric space in which
there is a unique arc connecting each pair of  points.
A group $\G$ is said to have 
{\textit{property $\FR$}} if for every action of $\G$ by  isometries
on an $\R$--tree, the fixed point set $\fix(\G)$
is nonempty. Serre's property $\FA$ is similar except that one
considers only actions on simplicial trees. A group
has $\FA$ if and only if it cannot be
expressed as a nontrivial amalgamated free product or
HNN extension.

\begin{lemmaa}\label{lem1} Let $\Gamma$ be a group 
that is generated by the union of the
subsets $A_1,\dots,A_N$.
If $H_{i,j}=\langle
A_i\cup A_j\rangle$ has property $\FR$ for 
all $i,j\in\{1,\dots,N\}$, then
$\Gamma$ also has property $\FR$.
\end{lemmaa}

\begin{proof} Let $C_1,\dots,C_N$ be connected subsets of an
$\R$--tree. It is not difficult to see that if
 $C_i\cap C_j$ is nonempty for $i,j=1,\dots,N$,
then $C_1\cap\dots\cap C_N$ is nonempty (cf Serre~\cite[p~65]{serre}).
Setting $C_i:=\fix(A_i)$ proves the lemma, since 
$C_i\cap C_j = \fix(H_{i,j})$ is assumed to be nonempty and
 $C_1\cap\dots\cap C_N = \fix (\G)$.
\end{proof}

Every finite group $G$ has $\FR$ because
the circumcentre of any $G$--orbit in an $\R$--tree will be a 
fixed point.

\begin{corolla}[The Triangle Criterion]
If $\Gamma$ is generated by $A_1\cup A_2\cup A_3$ 
and  $H_{i,j}=\langle
A_i\cup A_j\rangle$ is finite for $i,j=1,2,3$, then
$\Gamma$ has property $\FR$.
\end{corolla}

Let $\An$ denote the automorphism group
of the free group of rank $n$.

\begin{thma}\label{main} If $n\ge 3$ then $\An$
and $\slnz$  satisfy the Triangle Criterion and hence have property $\FR$.
\end{thma}

J-P\,Serre \cite{serre} was the first to prove that $\slnz$ has FA
if $n\ge3$,
and his argument shows that these groups actually have $\FR$.
Our argument is very similar to his except that he exploited
the pattern of nilpotent subgroups rather than finite ones. 
In the light of a theorem of J\,Tits \cite{tits}, Serre's argument
shows that all subgroups of finite index in 
$\slnz$ have FA. In contrast, there is a subgroup of finite index
$\G\subset\Aut(F_3)$ that does not have property $\FA$ 
(see McCool \cite{McC}), while it is 
unknown if $\Aut(F_n)$ has such subgroups when $n\ge 4$.
O\,Bogopolski \cite{Bog} was the first to prove that $\autn$
has FA.
M\,Culler and K\,Vogtmann \cite{CV} gave a short proof 
based on their idea of
``minipotent" elements.

The obvious appeal of \fullref{main} lies in the final phrase, 
but the stronger fact that these groups satisfy the Triangle Criterion
is useful in my work on fixed point theorems for
actions of automorphism groups of free groups on higher-dimensional
CAT$(0)$ spaces \cite{helly}. One can extend the theorem in
various ways (cf~\fullref{variation}) but I shall
not present the details here as to do so would 
obscure the simple and transparent proof that $\An$
has property FA, which is the main point of this note.
I hope that it is a proof that Zieschang would have
enjoyed. 

\subsubsection*{Acknowledgements} This research was supported
by Fellowships from the EPSRC and by a Royal
Society Wolfson Research Merit Award.

\section[Generating Aut(Fn) and SL(n,Z) by finite subgroups]{Generating $\autn$ and $\sln$ by finite subgroups}

We assume that $n\ge 3$
and  fix a basis $\B=\{a_1,\ldots,a_n\}$ of $F_n$. For
$i=1,\dots,n$, let $\e_i$ be the
automorphism of $F_n$ that sends $a_i$ to  $a_i^{-1}$
and fixes the other basis elements. J\,Nielsen proved that $\autn$ is generated
by the right  Nielsen
transformations $\rho_{ij}\co 
[a_i\mapsto a_ia_j,\ a_k\mapsto a_k
\text{ if } k\neq i]$ and the involutions $\e_i$. 

Let $\S_n\subset \autn$
be the group generated by permutations\footnote{We shall
write $(a_i\, a_j)$ to
denote the transposition of $a_i$ and $a_j$.} of $\B$. 
Conjugation by a permutation
$\sigma$ sends $\rho_{ij}$ to $\rho_{\sigma(i)\sigma(j)}$
and $\e_i$ to $\e_{\sigma(i)}$. Therefore 
$\autn$ is generated by 
$\rho_{12},\, \S_n$ and $\e_n$. In particular $\autn$ is generated
by $\rho_{12}\circ\e_2$ and the subgroup $W_n\cong (\Z_2)^n\rtimes \S_n$
generated by $\S_n$ and the $\e_i$. (The action of $\autn$ on the abelianisation
of $F_n$ gives a epimorphism $\autn \to \gln$, and the
image of $W_n$ under this map is the group of monomial matrices.)

We write
$\S_{n-2}\subset\S_n$ and $W_{n-2}\subset W_n$
for the subgroups
corresponding to the sub-basis $\{a_3,\dots,a_n\}$. Let
 $\theta :=\rho_{12}\circ\e_2$, let
$\tau := (a_2\, a_3)\circ\e_1$ and $\eta:=(a_1\, a_2)\circ\e_1\circ\e_2$,
 and note that each is an involution. Define
$$
A_1=\{\e_n, \eta\}\cup \S_{n-2},\ \ A_2=\{\theta\},\ \ A_3=\{\tau\}.
$$

\begin{lemma} $\autn$ is generated by $A_1\cup A_2\cup A_3$.
\end{lemma}

\begin{proof} Conjugating $(a_n\, a_3)\in\S_{n-2}$ by $\tau$
we get $(a_n\, a_2)$, which conjugates $\e_n$ to $\e_2$ and
$(a_n\, a_3)$ to $(a_2\,a_3)$. Thus $\e_1=(a_2\,a_3)\circ\tau$
and $(a_1\,a_2)=\eta\circ\e_1\circ\e_2$ are in
the subgroup generated
by the $A_i$; hence  $\S_n$ and $W_n$ are too. We already
noted that $\autn$ is generated by $W_n$ and $\theta$.
\end{proof}

\begin{lemma} The groups $H_{ij}=\langle A_i\cup A_j\rangle$ are finite.
\end{lemma}

\begin{proof}
$\S_{n-2}$ and $\e_n$
commute with the involutions $\theta$ and $\eta$,
and $\theta\circ\eta$ has order $3$, so
$H_{12}\cong W_{n-2}\times D_6$. As
$(\theta\circ\tau)^4=1$, we have $H_{23}\cong D_8$. And
$H_{13}\subset W_n$.
\end{proof}

These lemmas
prove that $\autn$ satisfies the Triangle Criterion if $n\ge 3$, and
by taking the images  of the $A_i$ under the natural map
$\autn\to \gln$ we see that $\gln$ does too. When $n$ is odd,
we obtain the corresponding result for $\sln$ by replacing the
image $\bar\gamma$ of each $\gamma\in A_i$ by
${\rm{det}}(\bar\gamma).\bar\gamma$; let $A_i^+(n)$
denote the image of $A_i$ modified in this manner.

When $n\ge 4$ is
even we need to adjust the $A_i$ a little more. 
Let $\alpha = \e_n\circ (a_n\, a_{n-1})$
and note that $\sln$ is generated
by the image $\bar\alpha$ of $\alpha$ and the subgroup
${\text{\rm{SL}}(n-1,\mathbb Z)}\subset \sln$ corresponding to the
sub-basis $\{a_1,\dots,a_{n-1}\}$.  
If $n\ge 4$ then the groups $H_{12}$ and $H_{13}$ remain
finite if we add $\alpha$ to $A_1$. Thus the sets $A_1^+(n-1)\cup
\{\alpha\}$, \mbox{$A_2^+(n-1)$}, $A_3^+(n-1)$ demonstrate that
$\sln$ satisfies the
Triangle Criterion.

\begin{rem} If $n\ge 4$ then by modifying the
sets $A_i$ slightly one can also show that ${\rm{SAut}}(F_n)$, the inverse image in
$\autn$ of $\sln$, satisfies the Triangle Criterion. 
\end{rem}

\subsection[The geometry of the Hij]{The geometry of the $H_{ij}$}

It would be unfair of me to leave the reader to guess
the origin of the finite subgroups used in the above proof,
so let me explain the geometry behind the construction.

Any finite subgroup of $\autn$ can be realised as a group of 
basepoint-preserving
isometries of a graph of Euler characteristic $1-n$. 
\fullref{figure}\footnote{I am grateful to Karen
Vogtmann for producing this figure.}
 below gives such realisations $Y_{ij}$ for the groups
$H_{ij}$. An important point is that if $i\notin\{j,k\}$ then
$A_i$ cannot be realised as a group of symmetries of $Y_{jk}$.
I wanted to obtain the generating set $W_n\cup\{\theta\}$
that proved useful in my work with Karen Vogtmann \cite{BV}. Thus,
starting with the rose and the graph $Y_{12}$ for $\theta$, I
looked for a third graph where $\theta$ could be realised together
with a symmetry intertwining $\{a_1,a_2\}$ and $\{a_3,\dots,a_n\}$.

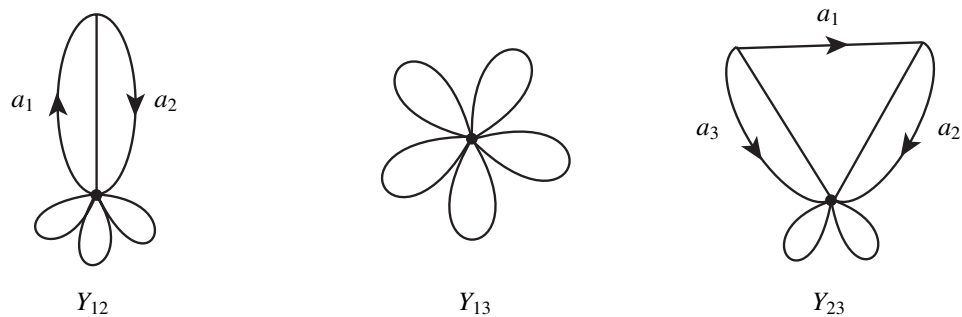
\begin{figure}[ht!]
\labellist\small
\pinlabel $a_1$ [r] at 11 111
\pinlabel $a_2$ [l] at 55 111
\pinlabel $Y_{12}$ [t] at 32 33
\pinlabel $Y_{13}$ [t] at 195 33
\pinlabel $Y_{23}$ [t] at 346 33
\pinlabel $a_3$ [r] at 303 99
\pinlabel $a_1$ [b] at 346 141
\pinlabel $a_2$ [l] at 389 99
\endlabellist
\centering
\includegraphics[width=.95\textwidth]{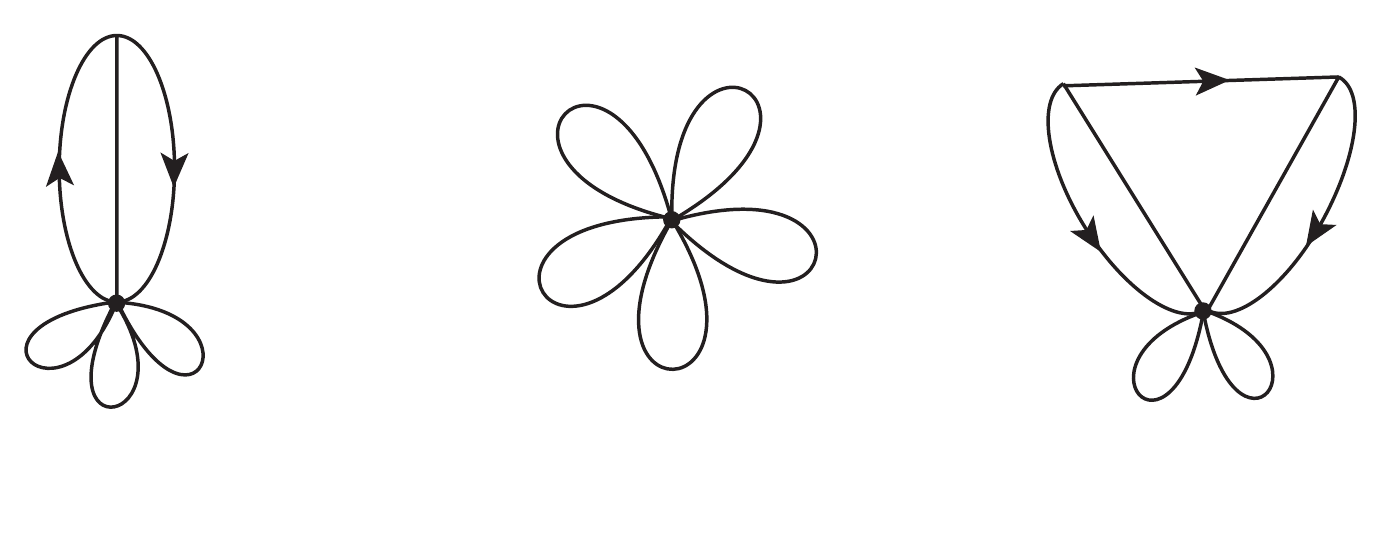}
\caption{The graphs $Y_{ij}$ exhibiting the finiteness of $H_{ij}$}
\label{figure}
\end{figure}

\section{Variations on the theme}\label{variation}

I have concentrated on configurations of finite subgroups in
this note but \fullref{lem1} can also be applied to situations
where the subgroups $\langle A_i\rangle$ are infinite. For example, if
$\gamma\in\Gamma$ lies in the
commutator subgroup of its centralizer, then $\fix(\gamma)$ 
will be nonempty whenever $\Gamma$ acts by isometries on an
$\mathbb R$--tree. (This is a special instance of a general
fact about semisimple actions on CAT$(0)$ spaces \cite{BH}.)
By exploiting such facts in conjunction with \fullref{lem1}
one can prove, for example, that the mapping class group of a
surface of genus at least 3 has property $\FR$, a result first
proved in \cite{CV}.

One can also strengthen \fullref{lem1} using an argument due
to J-P\,Serre \cite[p~64]{serre}: it suffices to require that
the $A_i$ have $\FR$ and, in any action of $\G$ on an $\R$--tree,  
that $a_ia_j$ have a fixed point, for every $a_i\in A_i$ 
and $a_j\in A_j$. To see this, one reduces to the case
 $n=2$ and argues that
if the fixed point sets of $A_1$ and $A_2$ did not 
intersect then  the point of $\fix(A_1)$ closest to $\fix(A_2)$ would be
fixed by all $a_1a_2$ with $a_1\in A_1$ and $a_2\in A_2$, which is 
a contradiction.

A quite different strengthening begins with the observation that
the behaviour of convex sets described in the proof of
\fullref{lem1} is a manifestation of the fact
that trees are $1$--dimensional objects. A suitable version 
of Helly's Theorem provides constraints on the
way in which convex sets can intersect in higher-dimensional
CAT$(0)$ spaces, and by applying these constraints to the
fixed point sets of finite subgroups one can prove far-reaching
generalisations of \fullref{main}; this is the theme of \cite{helly}.

\bibliographystyle{gtart}
\bibliography{link}

\end{document}